# Modeling inequality and spread in multiple regression[*]

## Rolf Aaberge[1], Steinar Bjerve[2] and Kjell Doksum[3]

*Statistics Norway, University of Oslo and University of Wisconsin, Madison*

**Abstract:** We consider concepts and models for measuring inequality in the distribution of resources with a focus on how inequality varies as a function of covariates. Lorenz introduced a device for measuring inequality in the distribution of income that indicates how much the incomes below the $u^{th}$ quantile fall short of the egalitarian situation where everyone has the same income. Gini introduced a summary measure of inequality that is the average over $u$ of the difference between the Lorenz curve and its values in the egalitarian case. More generally, measures of inequality are useful for other response variables in addition to income, e.g. wealth, sales, dividends, taxes, market share and test scores. In this paper we show that a generalized van Zwet type dispersion ordering for distributions of positive random variables induces an ordering on the Lorenz curve, the Gini coefficient and other measures of inequality. We use this result and distributional orderings based on transformations of distributions to motivate parametric and semiparametric models whose regression coefficients measure effects of covariates on inequality. In particular, we extend a parametric Pareto regression model to a flexible semiparametric regression model and give partial likelihood estimates of the regression coefficients and a baseline distribution that can be used to construct estimates of the various conditional measures of inequality.

## 1. Introduction

Measures of inequality provide quantifications of how much the distribution of a resource $Y$ deviates from the egalitarian situation where everyone has the same amount of the resource. The coefficients in location or location-scale regression models are not particularly informative when attention is turned to the influence of covariates on inequality. In this paper we consider regression models that are not location-scale regression models and whose coefficients are associated with the effect of covariates on inequality in the distribution of the response $Y$.

We start in Section 2.1 by discussing some familiar and some new measures of inequality. Then in Section 2.2 we relate the properties of these measures to a statistical ordering of distributions based on transformations of random variables that is equivalent to defining the distribution $H$ of the response $Z$ to have more resource inequality than the distribution $F$ of $Y$ if $Z$ has the same distribution as $q(Y)Y$ for some positive nondecreasing function $q(\cdot)$. Then we show that this ordering implies the corresponding ordering of each measure of inequality. We also consider orderings of distributions based on transformations of distribution functions and relate them to inequality. These notions and results assist in the construction of regression models with coefficients that relate to the concept of inequality.

---

[*]We would like to thank Anne Skoglund for typing and editing the paper and Javier Rojo and an anonymous referee for helpful comments. Rolf Aaberge gratefully acknowledges ICER in Torino for financial support and excellent working conditions and Steinar Bjerve for the support of The Wessmann Society during the course of this work. Kjell Doksum was supported in part by NSF grants DMS-9971301 and DMS-0505651.

[1]Research Department, Statistics Norway, P.O. Box 813, Dep., N-0033, Oslo, Norway, e-mail: `Rolf.Aaberge@ssb.no`

[2]Department of Mathematics, University of Oslo, P.O. Box 1053, Blindern, 0316, Oslo, Norway, e-mail: `steinar@math.uio.no`

[3]Department of Statistics, University of Wisconsin, 1300 University Ave, Madison, WI 53706, USA, e-mail: `doksum@stat.wisc.edu`

Section 3 shows that scaled power transformation models with the power parameter depending on covariates provide regression models where the coefficients relate to the concept of resource inequality. Two interesting particular cases are the Pareto and the log normal transformation regression models. For these models the Lorenz curve for the conditional distribution of Y given covariate values takes a particularly simple and intuitive form. We discuss likelihood methods for the statistical analysis of these models.

Finally, in Section 4 we consider semiparametric Lehmann and Cox type models that are based on power transformations of a baseline distribution $F_0$, or of $1 - F_0$, where the power parameter is a function of the covariates. In particular, we consider a power transformation model of the form

$$(1.1) \qquad F(y) = 1 - (1 - F_0(y))^{\alpha(\boldsymbol{x})},$$

where $\alpha(\boldsymbol{x})$ is a parametric function depending on a vector $\boldsymbol{\beta}$ of regression coefficients and an observed vector of covariates $\boldsymbol{x}$. This is an extension of the Pareto regression model to a flexible semiparametric model. For this model we present theoretical and empirical formulas for inequality measures and point out that computations can be based on available software.

## 2. Measures of inequality and spread

### 2.1. Defining curves and measures of inequality

The Lorenz curve (LC) is defined (Lorenz [19]) to be the proportion of the total amount of wealth that is owned by the "poorest" $100 \times u$ percent of the population. More precisely, let the random income $Y > 0$ have the distribution function $F(y)$, let $F^{-1}(y) = \inf\{y : F(y) \geq u\}$ denote the left inverse, and assume that $0 < \mu < \infty$, where $\mu = E(Y)$. Then the $LC$ (see e.g. Gastwirth [14]) is defined by

$$(2.1) \qquad L(u) = L_F(u) = \mu^{-1} \int_0^u F^{-1}(s)ds, \quad 0 \leq u \leq 1.$$

Let $I_{\{A\}}$ denote the indicator of the event $A$. For $F$ continuous we can write

$$(2.2) \qquad L(u) = \mu^{-1} E\{Y I_{\{Y \leq F^{-1}(u)\}}\}.$$

When the population consists of incomes of people, the $LC$ measures deviation from the *egalitarian* case $L(u) = u$ corresponding to where everyone has the same income $a > 0$ and the distribution of $Y$ is degenerate at $a$. The other extreme occurs when one person has all the income which corresponds to $L(u) = 0$, $0 \leq u \leq 1$. The intermediate case where $Y$ is uniform on $[0, b], b > 0$, corresponds to $L(u) = u^2$. In general $L(u)$ is non-decreasing, convex, below the line $L(u) = u$, $0 \leq u \leq 1$, and the greater the "distance" from $u$,ated the greater is the inequality in the population. If the population consists of companies providing a certain service or product, the $LC$ measures to what extent a few companies dominate the market with the extreme case corresponding to monopoly.

A closely related curve is the Bonferroni curve ($BC$) $B(u)$ which is defined (Aaberge [1], [2]), Giorgi and Mondani [15], Csörgö, Gastwirth and Zitikis [11]) as

$$(2.3) \qquad B(u) = B_F(u) = u^{-1}L(u),\ 0 \leq u \leq 1.$$

When $F$ is continuous the $BC$ is the $LC$ except that truncation is replaced by conditioning

$$(2.4) \qquad B(u) = \mu^{-1} E\{Y | Y \leq F^{-1}(u)\}.$$

The $BC$ possesses several attractive properties. First, it provides a convenient alternative interpretation of the information content of the Lorenz curve. For a fixed $u$, $B(u)$ is the ratio of the mean income of the poorest $100 \times u$ percent of the population to the overall mean. Thus, the $BC$ may also yield essential information on poverty provided that we know the poverty rate. Second, the $BC$ of a uniform $(0,a)$ distribution proves to be the diagonal line joining the points $(0,0)$ and $(1,1)$ and thus represents a useful reference line, in addition to the two well-known standard reference lines. The egalitarian reference line coincides with the horizontal line joining the points $(0,1)$ and $(1,1)$. At the other extreme, when one person holds all income, the $BC$ coincides with the horizontal axis except for $u = 1$.

In the next subsection we will consider ordering concepts from the statistics literature. Those concepts motivate the introduction of the following measures of concentration

$$(2.5) \qquad C(u) = C_F(u) = \int_0^u \left[\frac{F^{-1}(s)}{F^{-1}(u)}\right] ds = \mu_F \frac{L_F(u)}{F^{-1}(u)},\ 0 < u < 1$$

and

$$(2.6) \qquad D(u) = D_F(u) = \frac{1}{u} \int_0^u \left[\frac{F^{-1}(s)}{F^{-1}(u)}\right] ds = \mu_F \frac{B_F(u)}{F^{-1}(u)},\ 0 < u < 1.$$

Accordingly, *D(u)* emerges by replacing the overall mean $\mu$ in the dominator of *B(u)* by the $u^{th}$ quantile $y_u = F^{-1}(u)$ and is equal to the ratio between the mean income of those with lower income than the $u^{th}$ quantile and the $u$-quantile income. Thus, *C(u)* and *D(u)* measure inequality in income below the $u^{th}$ quantile. They satisfy $C(u) \leq u$, $D(u) \leq 1$, $0 < u < 1$, and $C(u)$ equals $u$ and $0$ while $D(u)$ equals $1$ and $0$ in the egalitarian and extreme non-egalitarian cases, respectively, and they equal $u/2$ and $1/2$ in the uniform case.

To summarize the information content of the inequality curves we recall the following inequality indices

$$(2.7) \qquad G = 2\int_0^1 \{u - L(u)\}du\ \text{(Gini)},\ B = \int_0^1 \{1 - B(u)\}du\ \text{(Bonferroni)},$$

$$(2.8) \qquad C = 2\int_0^1 \{u - C(u)\}du,\ D = \int_0^1 \{1 - D(u)\}du.$$

These indices measure distances from the curves to their values in the egalitarian case, take values between 0 and 1 and are increasing with increasing inequality. If all units have the same income then $G = B = C = D = 0$, and in the extreme non-egalitarian case where one unit has all the income and the others zero, $G = B = C = D = 1$. When $F$ is uniform on $[0, b]$, $B = C = D = 1/2$ and $G = 1/3$. The inequality curves $L(u), B(u), C(u), D(u)$, and *the inequality measures $G, B, C$ and $D$ are scale invariant*; that is, they remain the same if $Y$ is replaced by $aY, a > 0$.

### 2.2. Ordering inequality by transforming variables

When we are interested in how covariates influence inequality we may ask whether larger values of a covariate lead to more or less inequality. For instance, is there less inequality among the higher educated? To answer such questions we consider orderings of distributions on the basis of inequality, see e.g. Atkinson [5], Shorrocks and Foster [26], Dardanoni and Lambert [12], Muliere and Scarsini [20], Yitzhaki and Olkin [29], Zoli [30], and Aaberge [3]. In statistics and reliability engineering, orderings are plentiful, e.g. Lehmann [18], van Zwet [27], Barlow and Prochan [6], Birnbaum, Esary and Marshall [9], Doksum [13], Yanagimoto and Sibuya [28], Bickel and Lehmann [7], [8], Rojo and He [21], Rojo [22] and Shaked and Shanthikumar [25]. In statistics, similar orderings are often discussed in terms of spread or dispersion. Thus, for non-negative random variables, we could define $Y$ to have a distribution which is more spread out to the right than that of $Y_0$ if $Y$ can be written as $Y = h(Y_0)$ for some non-negative, nondecreasing convex function $h$ (using van Zwet [27]). It turns out to be more general and more convenient to replace "convex" with "starshaped" (convex functions $h$ are starshaped and concave functions $g$ are anti-starshaped provided $g(0) = h(0) = 0$).

Recall that a nondecreasing function $g$ defined on the interval $I \subset [0, \infty)$, is *starshaped* on $I$ if $g(\lambda x) \leq \lambda g(x)$ whenever $x \in I, \lambda x \in I$ and $0 \leq \lambda \leq 1$. Thus if $I = (0, \infty)$, for any straight line through the origin, then the graph of $g$ initially lies on or below it, and then lies on or above it. If $g(\lambda x) \geq \lambda g(x)$, $g$ is *anti-starshaped*. On the class $\mathcal{F}$ of continuous and strictly increasing distributions $F$ with $F(0) = 0$, Doksum [13] introduced the following partial ordering $F <_* H$ ($F$ is *starshaped with respect to $H$*) if $H^{-1}F$ is starshaped on $\{x : 0 < F(x) < 1\}$. This ordering was also considered by Yanagimoto and Sibuya [28] and Bickel and Lehmann [7,8]. Note that if $F <_* H$ and $X$ has distribution $F$, then $Z = H^{-1}[F(X)]$ has distribution $H$ and is a starshaped transformation of X. Moreover,

**Proposition 2.1.** *Suppose that $X$ and $Z$ have distributions $F$ and $H$, where $F, H \in \mathcal{F}$. Then $F <_* H$ if and only if there exists a positive nondecreasing function $q(\cdot)$ on $\{x : 0 < F(x) < 1\}$ such that $Z$ has the same distribution as $q(X)X$.*

*Proof.* Suppose that $F <_* H$; then $q(x) = H^{-1}(F(x))/x$ will do because the starshaped condition $g(\lambda x) \leq \lambda x$ is equivalent to $[g(\lambda x)/\lambda x] \leq g(x)/x$. That is,

(2.9) $$F <_* H \Leftrightarrow H^{-1}(F(x))/x \text{ is nondecreasing.}$$

Next, suppose that $q(\cdot)$ is positive and nondecreasing and that $Z = q(X)X$. Set $h(x) = q(x)x$. Then $h(x)$ is increasing and $P(X \leq x) = P(h(X) \leq h(x))$. It follows that $F(x) = H(q(x)x)$. That is, $q(x) = H^{-1}(F(x))/x$. □

Because of this proposition we say that if $F <_* H$ then $F$ is a *more egalitarian* distribution of resources than $H$.

We next show that the preceding definition of inequality leads to the corresponding ordering of the inequality curves $C_F(\cdot)$ and $D_F(\cdot)$ as well as of the indices $C$ and $D$.

**Proposition 2.2.** *Suppose that $F, H \in \mathcal{F}$ and $F <_* H$. Then $C_F(u) \geq C_H(u)$ and $D_F(u) \geq D_H(u)$, $0 < u < 1$. Therefore $C_F \leq C_H$ and $D_F \leq D_H$.*

*Proof.* It follows from (2.9), setting $u = F(x), v = F(x'), x < x'$, that $H^{-1}(u)/H^{-1}(v) \leq F^{-1}(u)/F^{-1}(v)$, for $0 < u < v < 1$. If we integrate this inequality over $u \in (0, v)$, we obtain $C_F(v) \geq C_H(v)$, $0 < v < 1$. The other inequalities follow from this. $\square$

The same order preservation as stated in Proposition 2.2 holds for $L(u)$ and $B(u)$.

**Theorem 2.1.** *Suppose that $F, H \in \mathcal{F}$ and $F <_* H$. Then $L_F(u) \geq L_H(u)$ and $B_F(u) \geq B_H(u)$, $0 < u < 1$. Moreover, $B_F \leq B_H$ and $G_F \leq G_H$.*

*Proof.* Let

$$\tag{2.10} a = \int_0^1 H^{-1}(v)dv \Big/ \int_0^1 F^{-1}(v)dv,$$

and consider the line $y = ax$ through the origin. Then $H^{-1}(F(x))$ initially lies on or below this line, say up to the point $x = b$. Thus

$$\tag{2.11} \begin{aligned} \int_0^u H^{-1}(v)dv &= \int_0^{F^{-1}(u)} H^{-1}(F(x))dF(x) \leq \int_0^{F^{-1}(u)} axdF(x) \\ &= a\int_0^u F^{-1}(v)dv \end{aligned}$$

for all $u \leq F^{-1}(b)$ which establishes $L_F(u) \geq L_H(u)$ for $u \leq F(b)$. On the other hand, for $x \geq b$, $y = H^{-1}(F(x))$ is above $y = ax$. Thus, for $u > F(b)$

$$s(u) = \int_b^{F^{-1}(u)} \left[ax - H^{-1}(F(x))\right]dF(x)$$

is a negative, decreasing function of $u$. We can write, for $u > F(b)$,

$$a\int_0^u F^{-1}(v)dv - \int_0^u H^{-1}(v)dv = a\int_0^{F(b)} F^{-1}(v)dv - \int_0^{F(b)} H^{-1}(v)dv + s(u)$$
$$\equiv c + s(u)$$

where $c$ is nonnegative by (2.11). It follows that $c + s(u)$ is a decreasing function that equals 0 when $u = 1$ by the definition of $a$. Thus, $c + s(u) \geq 0$ which establishes $L_F(u) \geq L_H(u)$ again by the definition of $a$. The other inequalities follow from this. $\square$

### 2.3. Ordering inequality by transforming distributions

A partial ordering on $\mathcal{F}$ based on transforming distributions rather than random variables is the following: $F$ represents *more equality than* $H$ ($F >_e H$) if

$$H(z) = g(F(z))$$

for some nonnegative increasing concave function $g$ on $[0, 1]$ with $g(0) = 0$ and $g(1) = 1$. In other words, $F <_e H$ if $F(H^{-1})$ is convex. If $F$ is uniform, the orderings $F >_e H$ implies $F <_* H$ and in this case the results of Propositions 2.1, 2.2 and Theorem 2.1 hold. Note that when $F$ and $H$ have densities $f$ and $h$, then $h(z) = g'(F(z))f(z)$ where $g'F$ is decreasing. That is, $F >_e H$ means that $F$ has relatively more probability mass on the right than $H$.

A similar ordering involves $\bar{F}(x) = 1 - F(x)$ and $\bar{H}(z) = 1 - H(z)$. In this case we say that $F$ represents a more equal distribution of resources than $H$ ($F >_r H$) if
$$\bar{H}(x) = g(\bar{F}(x))$$
for some nonegative increasing convex transformation $g$ on $[0, 1]$ with $g(0) = 0$ and $g(1) = 1$. In this case, if densities exist, they satisfy $h(z) = g'(\bar{F}(z))f(z)$, where $g'\bar{F}$ is decreasing. That is, relative to $F$, $H$ has mass shifted to the left.

**Remark.** Orderings of inequality based on transforming distributions can be restated in terms of orderings based on transforming random variables. Thus $F >_e H$ is equivalent to the distribution function of $V = F(Z)$ being convex when $X \sim F$ and $Z \sim H$.

## 3. Regression inequality models

### 3.1. Notation and introduction

Next consider the case where the distribution of $Y$ depends on covariates such as education, work experience, status of parents, sex, etc. Let $X_1, \ldots, X_d$ denote the covariates. We include an intercept term in the regression models, which makes it convenient to write $\boldsymbol{X} = (1, X_1, \ldots, X_d)^T$. Let $F(y|\boldsymbol{x})$ denote the conditional distribution of $Y$ given $\boldsymbol{X} = \boldsymbol{x}$ and define the quantile regression function as the left inverse of this distribution function. The key quantity is

$$\mu(u|\boldsymbol{x}) \equiv \int_0^u F^{-1}(v|\boldsymbol{x})dv.$$

With this notation we can write the regression versions of the Lorenz curve, for $0 < u < 1$ as
$$L(u|\boldsymbol{x}) = \mu(u|\boldsymbol{x})/\mu(1|\boldsymbol{x}), \ B(u|\boldsymbol{x}) = L(u|\boldsymbol{x})/u.$$

Similarly, $C(u|\boldsymbol{x})$, $D(u|\boldsymbol{x})$ and the summary coefficients $G(\boldsymbol{x})$, $B(\boldsymbol{x})$, $C(\boldsymbol{x})$ and $D(\boldsymbol{x})$ are defined by replacing $F(y)$ by $F(y|\boldsymbol{x})$. Note that estimates of $F(y|\boldsymbol{x})$ and $\mu(y|\boldsymbol{x})$ provide estimates of the regression versions of the curves and measures of inequality. Thus, the rest of the paper discusses regression models for $F(y|\boldsymbol{x})$ and $\mu(y|\boldsymbol{x})$. Using the results of Section 2, these models are constructed so that the regression coefficients reflect relationships between covariates and measures of inequality.

### 3.2. Transformation regression models

Let $Y_0$ with distribution $F_0$ denote a baseline variable which corresponds to the case where the covariate vector $\boldsymbol{x}$ has no effect on the distribution of income. We assume that $Y$ has a conditional distribution $F(y|\boldsymbol{x})$ which depends on $\boldsymbol{x}$ through some real valued function $\Delta(\boldsymbol{x}) = g(\boldsymbol{x}, \boldsymbol{\beta})$ which is known up to a vector $\boldsymbol{\beta}$ of unknown parameters. Let $Y \sim Z$ denote "$Y$ is distributed as $Z$". As we have seen in Section 2.2, if large values of $\Delta(\boldsymbol{x})$ correspond to a more egalitarian distribution of income than $F_0$, then it is reasonable to model this as

$$Y \sim h(Y_0),$$

for some increasing anti-starshaped function $h$ depending on $\Delta(\boldsymbol{x})$. On the other hand, an increasing starshaped $h$ would correspond to income being less egalitarian.

A convenient parametric form of $h$ is

(3.1) $$Y \sim \tau Y_0^\Delta,$$

where $\Delta = \Delta(\boldsymbol{x}) > 0$, and $\tau > 0$ does not depend on $\boldsymbol{x}$. Since $h(y) = y^{\Delta(\boldsymbol{x})}$ is concave for $0 < \Delta(\boldsymbol{x}) \le 1$, while convex for $\Delta(\boldsymbol{x}) > 1$, the model (3.1) with $0 < \Delta(\boldsymbol{x}) \le 1$ corresponds to covariates that lead to a less unequal distribution of income for $Y$ than for $Y_0$, while $\Delta(\boldsymbol{x}) \ge 1$ is the opposite case. Thus it follows from the results of Section 2.2 that if we use the parametrization $\Delta(\boldsymbol{x}) = exp(\boldsymbol{x}^T\boldsymbol{\beta})$, then the coefficient $\beta_j$ in $\boldsymbol{\beta}$ measures how the covariate $x_j$ relates to inequality in the distribution of resources $Y$.

**Example 3.1.** Suppose that $Y_0 \sim F_0$ where $F_0$ is the Pareto distribution $F_0(y) = 1 - (c/y)^a$, with $a > 1, c > 0, y \ge c$. Then $Y = \tau Y_0^\Delta$ has the Pareto distribution

(3.2) $$F(y|\boldsymbol{x}) = F_0\left(\left(\frac{y}{\tau}\right)^{\frac{1}{\Delta}}\right) = 1 - \left(\frac{\lambda}{y}\right)^{\alpha(\boldsymbol{x})}, \quad y \ge \lambda,$$

where $\lambda = c\tau$ and $\alpha(\boldsymbol{x}) = a/\Delta(\boldsymbol{x})$. In this case $\mu(u|\boldsymbol{x})$ and the regression summary measures have simple expressions, in particular

$$L(u|\mathbf{x}) = 1 - (1-u)^{1-\Delta(\mathbf{x})}.$$

When $\Delta(\boldsymbol{x}) = \exp(\boldsymbol{x}^T\boldsymbol{\beta})$ then $\log Y$ already has a scale parameter and we set $\alpha = 1$ without loss of generality. One strategy for estimating $\boldsymbol{\beta}$ is to temporarily assume that $\lambda$ is known and to use the maximum likelihood estimate $\hat{\boldsymbol{\beta}}(\lambda)$ based on the distribution of $\log Y_1, \ldots, \log Y_n$. Next, in the case where $(Y_1, X_1), \ldots, (Y_n, X_n)$ are i.i.d., we can use $\hat{\lambda} = n \min\{Y_i\}/(n+1)$ to estimate $\lambda$. Because $\hat{\lambda}$ converges to $\lambda$ at a faster than $\sqrt{n}$ rate, $\hat{\boldsymbol{\beta}}(\hat{\lambda})$ is consistent and $\sqrt{n}(\hat{\boldsymbol{\beta}}(\hat{\lambda}) - \boldsymbol{\beta})$ is asymptotically normal with the covariance matrix being the inverse of the $\lambda$-known information matrix.

**Example 3.2.** Another interesting case is obtained by setting $F_0$ equal to the log normal distribution $\Phi\big((\log(y) - \mu_0)/\sigma_0\big)$, $y > 0$. For the scaled log normal transformation model we get by straightforward calculation the following explicit form for the *conditional Lorenz curve*:

(3.3) $$L(u|\boldsymbol{x}) = \Phi\big(\Phi^{-1}(u) - \sigma_0 \Delta(\boldsymbol{x})\big).$$

In this case when we choose the parametrization $\Delta(\boldsymbol{x}) = \exp(x^T\boldsymbol{\beta})$, the model already includes the scale parameter $\exp(-\beta_0)$ for $\log Y$. Thus we set $\mu_0 = 1$. To estimate $\boldsymbol{\beta}$ for this model we set $Z_i = \log Y_i$. Then $Z_i$ has a $N\big(\alpha + \Delta(\boldsymbol{x}_i), \sigma_0^2\Delta^2(\boldsymbol{x}_i)\big)$ distribution, where $\alpha = \log \tau$ and $\mathbf{x}_i = (1, x_{i1}, \ldots, x_{id})^T$. Because $\sigma_0$ and $\alpha$ are unknown there are $d + 3$ parameters. When $Y_1, \ldots, Y_n$ are independent, this gives the log likelihood function (leaving out the constant term)

$$l(\alpha, \beta, \sigma_0^2) = -n\log(\sigma_0) - \sum_{i=1}^n \boldsymbol{x}_i^T\boldsymbol{\beta} - \frac{1}{2}\sigma_0^{-2}\sum_{i=1}^n \exp(-2\boldsymbol{x}_i^T\boldsymbol{\beta})\{Z_i - \alpha - \exp(\boldsymbol{x}_i^T\boldsymbol{\beta})\}^2$$

Likelihood methods will provide estimates, confidence intervals, tests and their properties. Software that only require the programming of the likelihood is available, e.g. Mathematica 5.2 and Stata 9.0.

## 4. Lehmann–Cox type semiparametric models. Partial likelihood

### *4.1. The distribution transformation model*

Let $Y_0 \sim F_0$ be a baseline income distribution and let $Y \sim F(y|\boldsymbol{x})$ denote the distribution of income for given covariate vector $\boldsymbol{x}$. In Section 2.3 it was found that one way to express that $F(y|\boldsymbol{x})$ corresponds to more equality than $F_0(y)$ is to use the model

$$F(y|\boldsymbol{x}) = h(F_0(y))$$

for some nonnegative increasing concave transformation $h$ depending on $\boldsymbol{x}$ with $h(0) = 0$ and $h(1) = 1$. Similarly, $h$ convex corresponds to a more egalitarian income. A model of the form $F_2(y) = h(F_1(y))$ was considered for the two-sample case by Lehmann [17] who noted that $F_2(y) = F_1^\Delta(y)$ for $\Delta > 0$ was a convenient choice of $h$. For regression experiments, we consider a regression version of this Lehmann model which we define as

(4.1) $$F(y|\boldsymbol{x}) = F_0^\Delta(y)$$

where $\Delta = \Delta(\boldsymbol{x}) = g(\boldsymbol{x},\boldsymbol{\beta})$ is a real valued parametric function and where $\Delta < 1$ or $\Delta > 1$ corresponds to $F(y|\boldsymbol{x})$ representing a more or less egalitarian distribution of resources than $F_0(y)$, respectively.

To find estimates of $\boldsymbol{\beta}$, note that if we set $U_i = 1 - F_0(Y_i)$, then $U_i$ has the distribution

$$H(u) = 1 - (1-u)^{\Delta(\boldsymbol{x})}, \; 0 < u < 1$$

which is the distribution of $F_0(Y_i)$ in the next subsection. Since the rank $R_i$ of $Y_i$ equals $N + 1 - S_i$, where $S_i$ is the rank of $1 - F_0(Y_i)$, we can use rank methods, or Cox partial likelihood methods, to estimate $\boldsymbol{\beta}$ without knowing $F_0$. In fact, because the Cox partial likelihood is a rank likelihood and $rank[1 - F_0(Y_i)] = rank(-Y_i)$, we can apply the likelihood in the next subsection to estimate the parameters in the current model provided we reverse the ordering of the $Y$'s.

### *4.2. The semiparametric generalized Pareto model*

In this section we show how the Pareto parametric regression model for income can be extended to a semiparametric model where the shape of the income distribution is completely general. This model coincides with the Cox proportional hazard model for which a wealth of theory and methods are available.

We defined a regression version of the Pareto model in Example 3.1 as

$$F(y|\boldsymbol{x}) = 1 - \left(\tfrac{c}{y}\right)^{\alpha_i}, \; y \geq c; \alpha_i > 0,$$

where $\alpha_i = \Delta_i^{-1}, \Delta_i = \exp\{\boldsymbol{x}_i^T \boldsymbol{\beta}\}$. This model satisfies

(4.2) $$1 - F(y|\boldsymbol{x}) = (1 - F_0(y))^{\alpha_i},$$

where $F_0(y) = 1 - c/y, y \geq c$. When $F_0$ is an arbitrary continuous distribution on $[0, \infty)$, the model (4.2) for the two sample case was called the Lehmann alternative by Savage [23], [24] because if $V$ satisfies model (4.1), then $Y = -V$ satisfies model (4.2). Cox [10] introduced proportional hazard models for regression experiments in survival analysis which also satisfy (4.2) and introduced partial likelihood methods that can be used to analyse such models even in the presence of censoring and time dependent covariates (in our case, wage dependent covariates).

Cox introduced the model equivalent to (4.2) as a generalization of the exponential model where $F_0(y) = 1 - \exp(-y)$ and $F(y|\boldsymbol{x_i})=F_0(\Delta_i^{-1}y)$. That is, (4.2) is in the Cox case a semiparametric generalization of a scale model with scale parameter $\Delta_i$. However, in our case we regard (4.2) as a semiparametric shape model which generalizes the Pareto model, and $\Delta_i$ represents the degree of inequality for a given covariate vector $\boldsymbol{x}_i$. The inequality measures correct for this confounding of shape and scale by being scale invariant.

Note from Section 2.3 that $\Delta_i < 1$ corresponds to $F(y|\boldsymbol{x})$ more egalitarian than $F_0(y)$ while $\Delta_i > 1$ corresponds to $F_0$ more egalitarian.

The Cox [10] partial likelihood to estimate $\boldsymbol{\beta}$ for (4.2) is (see also Kalbfleisch and Prentice [16], page 102),

$$L(\boldsymbol{\beta}) = \prod_{i=1}^{n} \left\{ \exp(-\boldsymbol{x}_{(i)}^T\boldsymbol{\beta}) \sum_{k \in R(Y_{(i)})} \exp(-\boldsymbol{x}_{(k)}^T\boldsymbol{\beta}) \right\}$$

where $Y_{(i)}$ is the $i$-th order statistic, $\boldsymbol{x}_{(i)}$ is the covariate vector for the subject with response $Y_{(i)}$, and $R(Y_{(i)}) = \{k : Y_{(k)} \geq Y_{(i)}\}$. Here $\hat{\boldsymbol{\beta}}=\arg\max L(\boldsymbol{\beta})$ can be found in many statistical packages, such as S-Plus, SAS, and STATA 9.0. These packages also give the standard errors of the $\hat{\beta}_j$. Note that $L(\boldsymbol{\beta})$ does not involve $F_0$.

Many estimates are available for $F_0$ in model (4.2) in the same packages. If we maximize the likelihood keeping $\boldsymbol{\beta} = \hat{\boldsymbol{\beta}}$ fixed, we find (e.g., Kalbfleisch and Prentice [16], p. 116, Andersen et al. [4], p. 483) $\hat{F}_0(Y_{(i)}) = 1 - \prod_{j=1}^{n} \hat{\alpha}_j$, where $\hat{\alpha}_j$ is the Breslow-Nelson-Aalen estimate,

$$\hat{\alpha}_j = \left\{ 1 - \frac{\exp(-\boldsymbol{x}_{(i)}^T\hat{\beta})}{\sum_{k \in R(Y_{(i)})} \exp(-\boldsymbol{x}_{(i)}^T\hat{\beta})} \right\}^{\exp(\boldsymbol{x}_{(i)}^T\beta)}$$

Andersen et al. [4] among others give the asymptotic properties of $\hat{F}_0$.

We can now give theoretical and empirical expressions for the conditional *inequality curves* and measures. Using (4.2), we find

(4.3) $$F^{-1}(u|\boldsymbol{x}_i) = F_0^{-1}(1 - (1-u)^{\Delta_i})$$

and

(4.4) $$\mu(u|\boldsymbol{x}_i) = \int_0^u F^{-1}(t|\boldsymbol{x}_i)dt = \int_0^u F_0^{-1}(1-(1-v)^{\Delta_i})dv.$$

We set $t = F_0^{-1}(1-(1-v)^{\Delta_i})$ and obtain

$$\mu(u|\boldsymbol{x}_i) = \Delta_i^{-1} \int_0^{\delta(u)} t(1-F_0(t))^{\Delta_i^{-1}-1}dF_0(t),$$

where $\delta_i(u) = F_0^{-1}(1 - (1-u)^{\Delta_i})$. To estimate $\mu(u|\boldsymbol{x_i})$, we let

$$b_i = \hat{F}_0(Y_{(i)}) - \hat{F}_0(Y_{(i-1)}) = \prod_{j=1}^{i-1} \hat{\alpha}_j = (1 - \hat{\alpha}_i)\prod_{j=1}^{i-1} \hat{\alpha}_j$$

be the jumps of $\hat{F}_0(\cdot)$; then

$$\hat{\mu}(u|\boldsymbol{x}_i) = \hat{\Delta}_i^{-1} \sum_j b_j Y_{(j)}(1 - \hat{F}_0(Y_{(j)}))^{\hat{\Delta}_i^{-1}-1}$$

where the sum is over $j$ with $\hat{F}_0(Y_{(j)}) \leq 1 - (1-u)^{\hat{\Delta}_i}$. Finally,

$$\hat{L}(u|\boldsymbol{x}) = \hat{\mu}(u|\boldsymbol{x})/\hat{\mu}(1|\boldsymbol{x}), \hat{B}(u|\boldsymbol{x}) = \hat{L}(u|\boldsymbol{x})/u,$$

and

$$\hat{C}(u|\boldsymbol{x}) = \hat{\mu}(u|\boldsymbol{x})/\hat{F}^{-1}(u|\boldsymbol{x}), \hat{D}(u|\boldsymbol{x}) = \hat{C}(u|\boldsymbol{x})/u,$$

where $\hat{F}^{-1}(u|\boldsymbol{x})$ is the estimate of the conditional quantile function obtained from (4.3) by replacing $\Delta_i$ with $\hat{\Delta}_i$ and $F_0$ with $\hat{F}_0$.

**Remark.** The methods outlined here for the Cox proportional hazard model have been extended to the case of ties among the responses $Y_i$, to censored data, and to time dependent covariates (see e.g. Cox [10], Andersen et al. [4] and Kalbfleisch and Prentice [16]). These extensions can be used in the analysis of the semiparametric generalized Pareto model with tied wages, censored wages, and dependent covariates.

## References

[1] AABERGE, R. (1982). On the problem of measuring inequality (in Norwegian). Rapporter 82/9, Statistics Norway.

[2] AABERGE, R. (2000a). Characterizations of Lorenz curves and income distributions. *Social Choice and Welfare* **17**, 639–653.

[3] AABERGE, R. (2000b). Ranking intersecting Lorenz Curves. Discussion Paper No. 412, Statistics Norway.

[4] ANDERSEN, P. K., BORGAN, Ø. GILL, R. D. AND KEIDING, N. (1993). *Statistical Models Based on Counting Processes.* Springer, New York. MR1198884

[5] ATKINSON, A. B. (1970). On the measurement of inequality, *J. Econ. Theory* **2**, 244–263. MR0449508

[6] BARLOW, R. E. AND PROSCHAN, F. (1965). *Mathematical Theory of Reliability.* Wiley, New York. MR0195566

[7] BICKEL, P. J. AND LEHMANN, E. L. (1976). Descriptive statistics for nonparametric models. III. Dispersion. *Ann. Statist.* **4**, 1139–1158. MR0474620

[8] BICKEL, P. J. AND LEHMANN, E. L. (1979). Descriptive measures for nonparametric models IV, Spread. In *Contributions to Statistics, Hajek Memorial Volume*, J. Juneckova (ed.). Reidel, London, 33–40. MR0561256

[9] BIRNBAUM, S. W., ESARY, J. D. AND MARSHALL, A. W. (1966). A stochastic characterization of wear-out for components and systems. *Ann. Math. Statist.* **37**, 816–826. MR0193727

[10] COX, D. R. (1972). Regression models and life tables (with discussion). *J. R. Stat. Soc. B* **34**, 187–220. MR0341758

[11] Csörgö, M., Gastwirth, J. L. and Zitikis, R. (1998). Asymptotic confidence bands for the Lorenz and Bonferroni curves based on the empirical Lorenz curve. *Journal of Statistical Planning and Inference* **74**, 65–91. MR1665121

[12] Dardanoni, V. and Lambert, P. J. (1988). Welfare rankings of income distributions: A role for the variance and some insights for tax reforms. *Soc. Choice Welfare* **5**, 1–17. MR0931516

[13] Doksum, K. A. (1969). Starshaped transformations and the power of rank tests. *Ann. Math. Statist.* **40**, 1167–1176. MR0243699

[14] Gastwirth, J. L. (1971). A general definition of the Lorenz curve. *Econometrica* **39**, 1037–1039.

[15] Giorgi, G. M. and Mondani, R. (1995). Sampling distribution of the Bonferroni inequality index from exponential population. *Sankyā* **57**, 10–18. MR1394865

[16] Kalbfleisch, J. D. and Prentice, R. L. (2002). *The Statistical Analysis of Failure Time Data*, 2nd edition. Wiley, New York. MR1924807

[17] Lehmann, E. L. (1953). The power of rank tests. *Ann. Math. Statist.* **24**, 23–43. MR0054208

[18] Lehmann, E. L. (1955). Ordered families of distributions. *Ann. Math. Statist.* **37**, 1137–1153. MR0071684

[19] Lorenz, M. C. (1905). Methods of measuring the concentration of wealth. *J. Amer. Statist.* **9**, 209–219.

[20] Muliere, P. and Scarsini, M. (1989). A Note on Stochastic Dominance and Inequality Measures. *Journal of Economic Theory* **49**, 314–323. MR1030876

[21] Rojo, J. and He, G. Z. (1991). New properties and characterizations of the dispersive orderings. *Statistics and Probability Letters* **11**, 365–372. MR1108359

[22] Rojo, J. (1992). A pure-tail ordering based on the ratio of the quantile functions. *Ann. Statist.* **20**, 570–579. MR1150363

[23] Savage, I. R. (1956). Contributions to the theory of rank order statistics – the two-sample case. *Ann. Math. Statist.* **27**, 590–615. MR0080416

[24] Savage, I. R. (1980). *Lehmann Alternatives*. Colloquia Mathematica Societatis János Bolyai, Nonparametric Statistical Inference, Proceedings, Budapest, Hungary.

[25] Shaked, M. and Shanthikumar, J. G. (1994). *Stochastic Orders and Their Applications*. Academic Press, San Diego. MR1278322

[26] Shorrocks, A. F. and Foster, J. E. (1987). Transfer sensitive inequality measures. *Rev. Econ. Stud.* **14**, 485–497. MR0896448

[27] van Zwet, W. R. (1964). *Convex Transformations of Random Variables*. Math. Centre, Amsterdam. MR0175265

[28] Yanagimoto, T. and Sibuya, M. (1976). Isotonic tests for spread and tail. *Annals of Statist. Math.* **28**, 329–342. MR0433711

[29] Yitzhaki, S. and Olkin, I. (1991). Concentration indices and concentration curves. Stochastic Order and Decision under Risk. *IMS Lecture Notes–Monograph Series.* MR1196066

[30] Zoli, C. (1999). Intersecting generalized Lorenz curves and the Gini index. *Soc. Choice Welfare* **16**, 183–196. MR1676972